\documentclass[onefignum,onetabnum]{siamonline220329}

\usepackage{lipsum}
\usepackage{amsfonts}
\usepackage{graphicx}
\usepackage{epstopdf}

\usepackage{enumitem}
\setlist[enumerate]{leftmargin=.5in}
\setlist[itemize]{leftmargin=.5in}

\newsiamremark{remark}{Remark}
\newsiamremark{hypothesis}{Hypothesis}
\crefname{hypothesis}{Hypothesis}{Hypotheses}
\newsiamthm{claim}{Claim}

\headers{Gradient Descent in Linear Regression}{K. Duraisamy}

\title{Finite Sample Analysis and Bounds of Generalization Error of Gradient Descent in In-Context Linear Regression\thanks{Submitted to the editors 05/03/2024.
\funding{This work was funded by AFOSR under contract FA9550-17-1-0195.}}}

\author{Karthik Duraisamy\thanks{University of Michigan, Ann Arbor, MI
  (\email{kdur@umich.edu}, \url{http://caslab.engin.umich.edu}).} }

\usepackage{amsopn}

\ifpdf
\hypersetup{
  pdftitle={Finite Sample Analysis of Generalization Error in Gradient Descent: Bounds and Error Rates in In-Context Linear Regression},
  pdfauthor={K. Duraisamy}
}
\fi
\usepackage{comment}
\usepackage{array}
\usepackage[margin=1in]{geometry}
\usepackage{amsmath}
\usepackage{amssymb}
\usepackage{hyperref}
\usepackage{subcaption}
\usepackage{tikz}
\usepackage{pgfplots}

\newcommand{\reals}{\mathbb{R}}

\DeclareMathOperator{\Tr}{\textrm{Tr}}
\newcommand{\norm}[1]{\| #1 \|_2}
\newcommand{\tr}{\textrm{Tr}}
\newcommand{\xh}{\hat{x}}

\usepackage{empheq}
\usepackage[most]{tcolorbox}

\newtcbox{\mymath}[1][]{%
	nobeforeafter, math upper, tcbox raise base,
	enhanced, colframe=blue!30!black,
	colback=blue!30, boxrule=1pt,
	#1}

\begin{document}
\maketitle

\begin{abstract}
Recent studies show that transformer-based architectures emulate gradient descent during a forward pass, contributing to in-context learning capabilities— an ability where the model adapts to new tasks based on a sequence of prompt examples without being explicitly trained or fine tuned to do so. This work investigates the generalization properties of a single step of gradient descent in the context of linear regression with well-specified models. A random design setting is considered and analytical expressions are derived for the statistical properties and bounds of generalization error in a non-asymptotic (finite sample) setting. These expressions are notable for avoiding arbitrary constants, and thus offer robust quantitative information and scaling relationships. These results are contrasted with those from classical least squares regression (for which analogous finite sample bounds are also derived), shedding light on systematic and noise components, as well as optimal step sizes. Additionally,  identities involving high-order products of Gaussian random matrices are presented as a byproduct of the analysis.
\end{abstract}



\section{Introduction}

Significant developments in large language models have led to much research  on enhancing and understanding  learning processes of Transformers~\cite{vaswani2017attention,khan2022transformers}. Among other characteristics, Transformer models have been shown to mimic the mechanisms of gradient descent during their forward pass~\cite{ICL1}. This behavior has been related to in-context learning abilities of Transformers~\cite{garg2022can}. In-context examples typically consist of input-output pairs that are directly related to a specific task. This mechanism allows the Transformer to process a new query input by leveraging the prompt examples to generate the corresponding predictive output effectively. The connections between in-context learning and gradient descent have been widely studied over the past two years. von Oswald et al.~\cite{ICL1} show that training Transformers on auto-regressive tasks mimics gradient-based meta-learning. In other words, Transformers learn in-context by emulating gradient descent, becoming meta-optimizers and excelling in regression tasks. Dai et al.~\cite{ICL4} hypothesize that attention values act as meta-gradients, enabling implicit fine-tuning for in-context learning. 

While the present work is strongly motivated by the above connections and hypotheses, the primary focus of this paper is on the foundational aspects of gradient descent method. We consider a case of in-context linear regression with a well-specified model, and aim to  investigate the extent to which a single step of gradient descent can generalize across examples from a noisy distribution. The focus on the single step is because of the relevance to in-context learning as described above.

The problem setting is as follows. Assume we have a model $y = W_0 x$, where for simplicity we assume $x \in \reals^n$ corresponds to realizations from the standard normal distribution. $W_0 \in \reals^{m \times n}$ is the existing weight matrix and $y \in \reals^m$ is the output.  As in the in-context setting, assume that we get N i.i.d. data pairs $\{x_1,y_1\},...\{x_N,y_N\}$ where $p(y_i|x_i) = \mathcal{N}(y_i;W_1 x_i,  \sigma^2 I)$.  We are interested in the properties of one-shot gradient descent, for an testing input $x_{N+1} \triangleq \hat{x}$. Particularly, we would like to know the statistical properties of the prediction $\hat{y}$ and bounds on the generalization error.  We will compare these results to least squares regression (on the given data pairs). A random design setting will be concerned (i.e. the in-context data available to us is drawn randomly).  

Indeed,  properties of gradient descent have been studied for a long time, with much of the work devoted to convergence properties (e.g. ~\cite{shapiro1996convergence,hu2016convergence}). Optimality of stochastic gradient descent has also been addressed from various perspectives (e.g. ~\cite{rakhlin2011making}). Ref.~\cite{neu2021information} presents  upper bounds on the generalization error that depend on local statistics of the stochastic gradients using  information-theoretic constructs. 
 We seek results that are non-asymptotic (i.e. by considering finite samples), not containing arbitrary constants, and not requiring bounded covariates. In fact, it can be argued that such results are hard to come by even in conventional linear least-squares regression. For instance,  Gyorfi et al.~\cite{gyorfi2002distribution} and Catoni's~\cite{catoni2004statistical} error bounds for least squares regression are well-crafted, but contain an arbitrary constant. Many publications contain terms such as $O(\cdot)$, thus rendering an asymptotic error estimate. Audibert \& Catoni~\cite{audibert2010linear} require 
$N >> n \log n$ and Hsu et al.~\cite{hsu2012random} require $N >> n$.  Classical PAC-Bayesian bounds~\cite{mcallester1999pac} require bounded loss functions or additional parameters beyond the data~\cite{germain2016pac}.  The author is careful to emphasize that the   above works are rigorous, and focused on a more general - and thus more practically relevant - setting than the present one, and that a bound can be useful even under the above conditions.

The outline and main contributions of the paper are as follows: The expected generalization error of in-context gradient descent (for linear regression with well-specified models) is derived in Section~\ref{sec:exgenerr}. Comparisons are made with classical least squares regression and a breakdown of the systematic and noise components and an expression is provided for the optimal step size .  In sections~\ref{sec:bounds} and ~\ref{sec:LSQ}, probabilistic bounds are derived for gradient descent and least squares regression. Section~\ref{sec:conclusions} explores connections to existing work. As a byproduct of this work, several identities were derived involving high order products of Gaussian random matrices, and provided in the Appendix.

\section{Expected Generalization Error}
\label{sec:exgenerr}
As mentioned in the introduction, the in-context learning setting is equated to one step of gradient descent over a `prompt' of $N$ i.i.d. data-pairs from a new task. Gradient descent  yields an output with the new weights, i.e. $\hat{y} = (W_0 + \Delta W)\hat{x}$.

\begin{theorem}[Expected Error]\label{thm:EE}
 Given $X \triangleq [x_1 \ \ x_2  \ \ ...\ \ x_N] \in \reals^{n \times N}$, which is a i.i.d standard normal random matrix and $Y \triangleq [y_1 \ \ y_2 \ \ .. \ \ y_N] \in \reals^{m \times N}$ with $y_i = W_1 x_i + \sigma^2 z_i$ and $z_i \in \reals^m$ is an i.i.d. sample from the standard normal distribution. For a step size $\eta$, the expected mean squared error for a step of gradient descent is 
$\mathbb{E}[\ell] =||W_1 - W_0 ||_F^2  \left((1-\eta)^2+ \eta^2 \frac{n+1}{N}\right) + \sigma^2 \left(m +  \eta^2 \frac{n}{N} \right)  $
\end{theorem} 
 
\begin{proof}For a mean squared training loss, it is easy to see that the update to the weights after one step of gradient descent is 

\begin{align*}
\Delta W &= - \frac{\eta}{N} \sum_i (W_0 x_i - y_i) x_i^T  \triangleq - \frac{\eta}{N} (W_0 X - Y) X^T. 
\end{align*}
Given $ y_i = W_1 x_i + \sigma z_i$,  
\begin{align*}
\Delta W &= - \frac{\eta}{N} (W_0 X - W_1 X { - \sigma Z}) X^T \\
 &= \frac{\eta}{N}(W_1-W_0) X X^T { + \sigma \frac{\eta}{N}Z  X^T},
\end{align*}
where $Z \triangleq [z_1 \ \ z_2 \ \ .. \ \ z_N] \in \reals^{m \times N}$. Given $X,Z$, it can be shown that the distribution of the weight update is $$\Delta W \sim \mathcal{W}\left(\eta(W_1-W_0),\frac{\eta^2}{N}(n+1)(W_1-W_0)^T(W_1-W_0)\right),$$
where $\mathcal{W}$ is the Wishart distribution.

We are interested in determining the statistical properties of the generalization error for the above random design case. We will begin with the expected mean squared error, which for a
given test location $\hat{x}$ is
$\ell  \triangleq \norm{W_1 \hat{x} + \sigma \hat{z} - (W_0 + \Delta W) \hat{x}}^2.$ Substituting $\Delta W$,
\begin{align*}
\ell &=  \norm{\left( W_1 \hat{x} - W_0 - \frac{\eta}{N}(W_1-W_0) X X^T \right)\hat{x}    { - \sigma \frac{\eta}{N}Z  X^T \hat{x}+ \sigma \hat{z}}}^2 \\
\mathbb{E}[\ell]&= \mathbb{E}[ \norm{\left((W_1-W_0)(I- \frac{\eta}{N} X X^T \right)\hat{x} }^2] {+ \sigma^2 \mathbb{E}[\hat{z}^T \hat{z}]   + \sigma^2 \frac{\eta^2}{N^2} \mathbb{E}\left[ \hat{x} X Z^T Z X^T \hat{x}\right]}\\
&= \mathbb{E}[ \norm{ A(I- a Q) \hat{x} }^2]  { + \sigma^2 m + \sigma^2 a^2 N n} \ \ \textrm{where} \ \ a \triangleq \frac{\eta}{N}, Q \triangleq X X^T, A \triangleq W_1-W_0 \\
\mathbb{E}_{\hat{x}}[\ell] &= \tr [(I- a Q)A^T A(I- a Q)] { + \sigma^2 m + \sigma^2 a^2 N n} \\
&= \tr [A^T A - a Q A^T A - A^T A a Q +a^2 Q A^T A Q] { + \sigma^2 m + \sigma^2 a^2 N n} \\
\mathbb{E}[\ell]&= \tr [A^T A](1 - a N - a N + a^2((n+1)N+N^2)) {+ \sigma^2 m + \sigma^2 a^2 N n} \\
&= \tr [A^T A](1 - 2 a N + a^2((n+1)N+N^2)) { + \sigma^2 m + \sigma^2 a^2 N n}.
\end{align*}
\end{proof}

\noindent Consider $m=1$  for simplicity. The above yields
\begin{equation}
\label{eq:ExpectedLoss}
\mathbb{E}[\ell]=  ||W_1 - W_0 ||_2^2  \left((1-\eta)^2+ \eta^2 \frac{n+1}{N}\right) { + \sigma^2 \left(1 +  \eta^2 \frac{n}{N} \right)}.
\end{equation}

Note that standard least squares regression yields (see Section~\ref{sec:LSQ})

$$
\mathbb{E}[\ell] = ||W_1-W_0||_2^2 \left( 1 - \frac{N}{n} \right) { + \sigma^2 \left( 1 + \frac{N}{n-N-1} \right)} \ \ \textrm{for}  \ \  2 \leq N \leq n-1 
$$
$$
\mathbb{E}[\ell] =  { \sigma^2 \left( 1 + \frac{n}{N-n-1} \right)} \ \ \textrm{for}  \ \  n+1 \leq N \le \infty. 
$$

In the under-parametrized regime (i.e. $N < n$), the systematic error in Least squares (or more precisely, Least Norm) regression drops linearly, whereas the error due to noise grows. The blow up near $N = n$ and the double descent behavior is attributable to noise amplification, and is - by now - well studied~\cite{belkin2019reconciling,belkin2020two}. For $N \ge n$, there is no systematic error as is well known.

\begin{figure}
	\centering
	\includegraphics[width=0.75\textwidth]{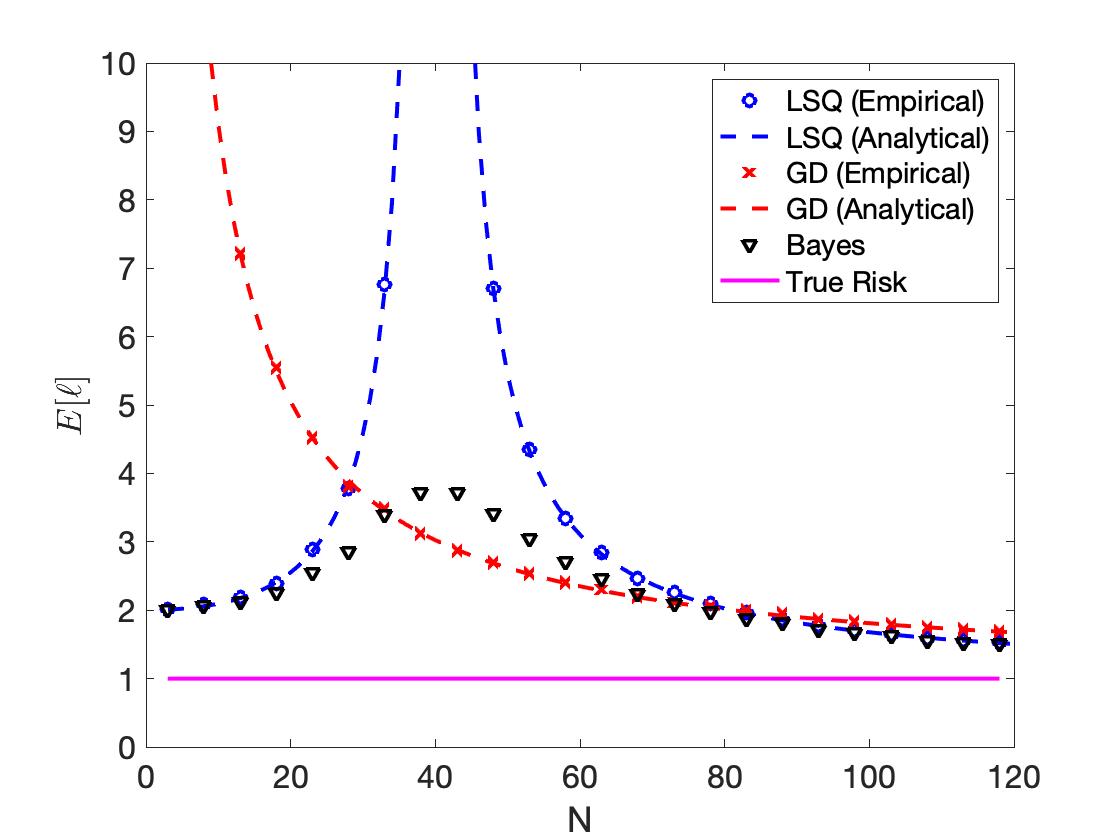}
		\includegraphics[width=0.49\textwidth]{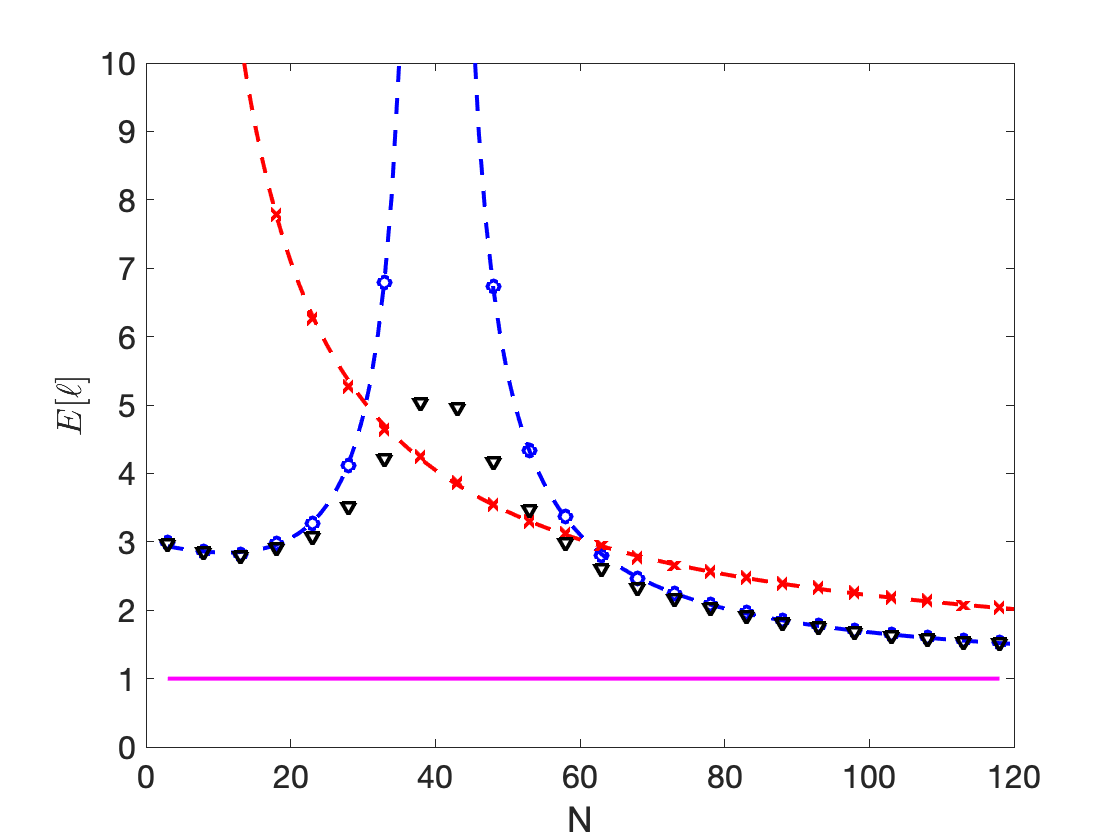}
				\includegraphics[width=0.49\textwidth]{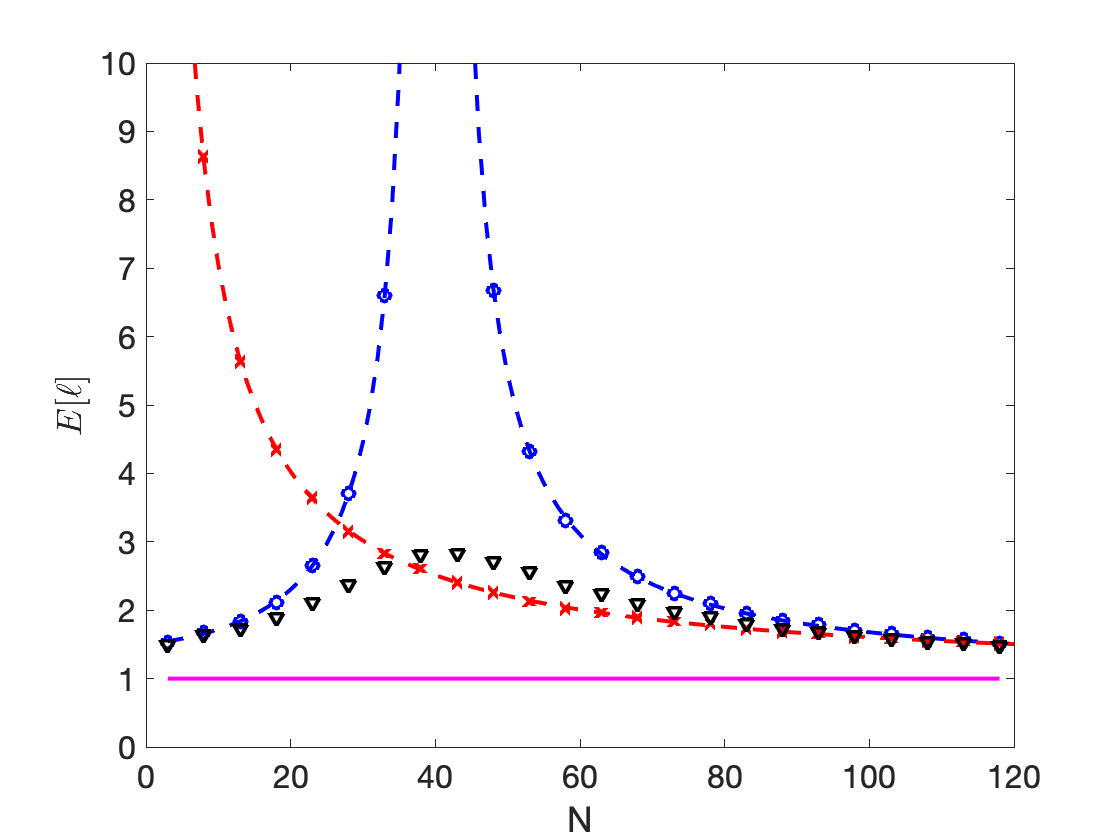}
	\caption{Expected Generalization Error for $n=40; m=1; \eta=1$. Top : $\norm{W_1}^2/\sigma^2=1$; Bottom left : $\norm{W_1}^2/\sigma^2=2$; Bottom right: $\norm{W_1}^2/\sigma^2=0.5$ }
	\label{fig:EL}
\end{figure}

 Figure ~\ref{fig:EL} shows for $n=40, m=1, \eta=1$, empirical evaluations (for each N, we used 500 random designs and 500 test evaluations) compared to the above analytical expressions. $W_1$ was sampled from the standard normal distribution and normalized such that $\norm{W_1} =1$ and $\sigma^2$ was varied to yield different signal-to-noise ratios.
 The breakdown of the different components of the testing error is also shown for $\norm{W_1}^2 = \sigma^2$.

\begin{figure}
	\centering
	\includegraphics[width=0.5\textwidth]{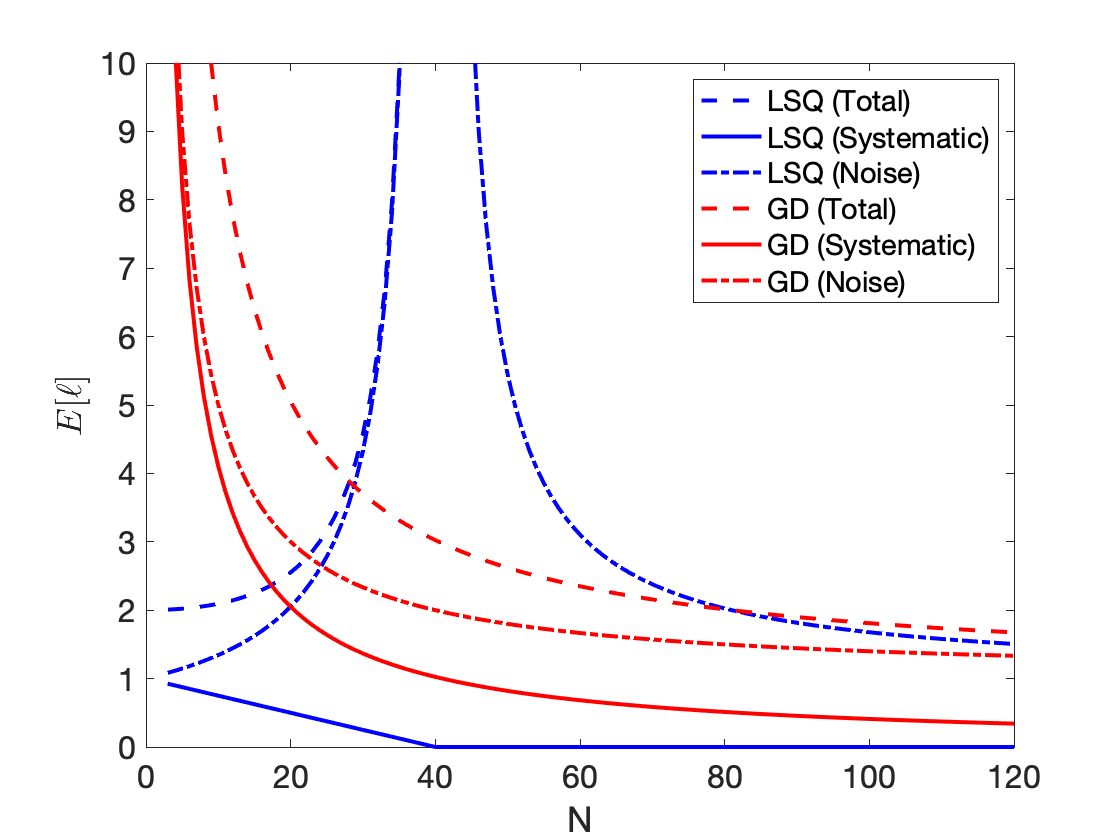}
	\caption{$n=40; m=1, \norm{W_1}^2/\sigma^2=1$}
	\label{fig:EL1}
\end{figure}

\begin{figure}
	\centering
	\includegraphics[width=0.5\textwidth]{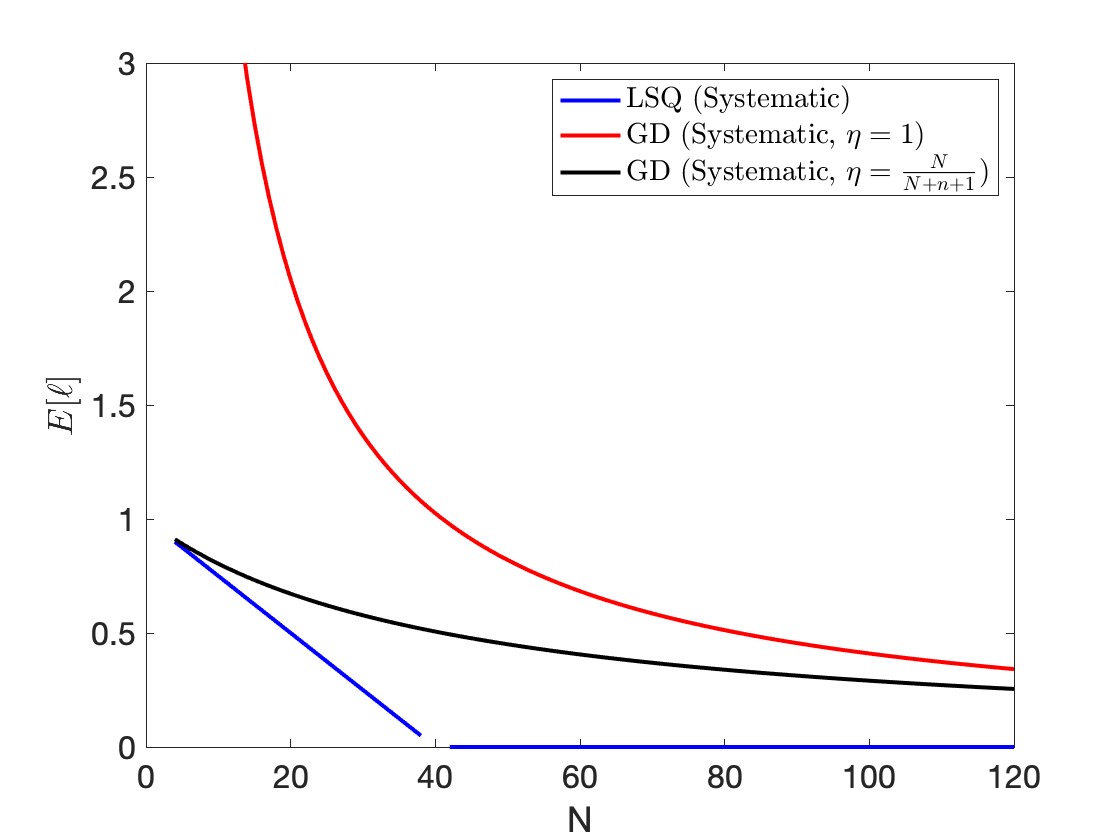}
	\caption{$n=40; m=1, \sigma=0$}
	\label{fig:ELopt}
\end{figure}

Examining Equation~\ref{eq:ExpectedLoss}, we can determine the optimal step size $$\eta_{opt} = \frac{N}{N+n+1 + \frac{\sigma^2}{||W_1||^2}n}$$ for gradient descent. In realistic scenarios, we do not  know the signal to noise ratio, and thus a practical guide would be 

$$
\eta_{opt} = \frac{N}{N+n+1} = 1 - \frac{n+1}{N+n+1}, 
$$
which for the noiseless case yields
$$\mathbb{E}[\ell] = ||W_1||^2 \left(1- \frac{N}{N+n+1}\right).$$ This is shown to greatly reduce the generalization error as shown in Fig~\ref{fig:ELopt}.

\section{Bounds}
\label{sec:bounds}
In the previous section, we examined the expected error in in-context regression. In this section, we derive probabilistic bounds for the generalization error.

\begin{theorem}[Bounds on Generalization Error]\label{thm:BB}
As in the prior setting, assume $X \triangleq $ $[x_1 \ \  ...\ \ x_N]$ $ \in \reals^{n \times N}$ as an i.i.d standard normal random matrix and $Y \triangleq [y_1 \ \ y_2 \ \ .. \ \ y_N] \in \reals^{1 \times N}$ with $y_i = w_1^T x_i + \sigma^2 z_i$ and $z_i \in \reals^m$ is an i.i.d. sample from the standard normal distribution.  With at least a probability of $1-\delta$,the generalization error of gradient descent (with a step size of 1) is bounded by

\begin{equation}
\ell \leq \left[\frac{||w_1-w_0||_2^2(n+1)  + \sigma^2 (n+N)}{N}  + \sqrt{\frac{V}{\delta}} \right], 
\end{equation}
where \begin{align*}
V & \triangleq  \sigma^4\left(2  + \frac{12 n+6 n^2}{N^3} + \frac{6 n + 2 n^2}{N^2} + \frac{4 n}{N}\right)\\
&+\sigma^2 \norm{w_1-w_0}^2\left(\frac{72  +  60 n  + 12 n^2}{N^3} + \frac{12 + 16 n + 4 n^2}{N^2} + \frac{4 + 4 n}{N}\right) \\
&+\norm{w_1-w_0}^4 \left( \frac{90  + 48 n + 6 n^2}{N^3} + \frac{20  + 10 n + 2 n^2}{N^2} \right)
\end{align*} 
\end{theorem}

\noindent Note: For the noiseless case, therefore with at least a probability of $1-\delta$, we can guarantee that the generalization error is bounded by
\begin{equation}
\ell \leq \left[\frac{n+1}{N}+ \sqrt{\frac{ \left[ \frac{90+48 n+6 n^2}{N^3} + \frac{20+ 10 n +2 n^2}{N^2}\right]}{\delta}} \right] \norm{w_1-w_0}^2.
\end{equation}

For large n, we have
\begin{equation}
\ell \leq \frac{n}{N}\left(1+ \sqrt{\frac{2}{\delta}} \right) \norm{w_1-w_0}^2.
\end{equation}

\begin{figure}
	\includegraphics[width=0.33\textwidth]{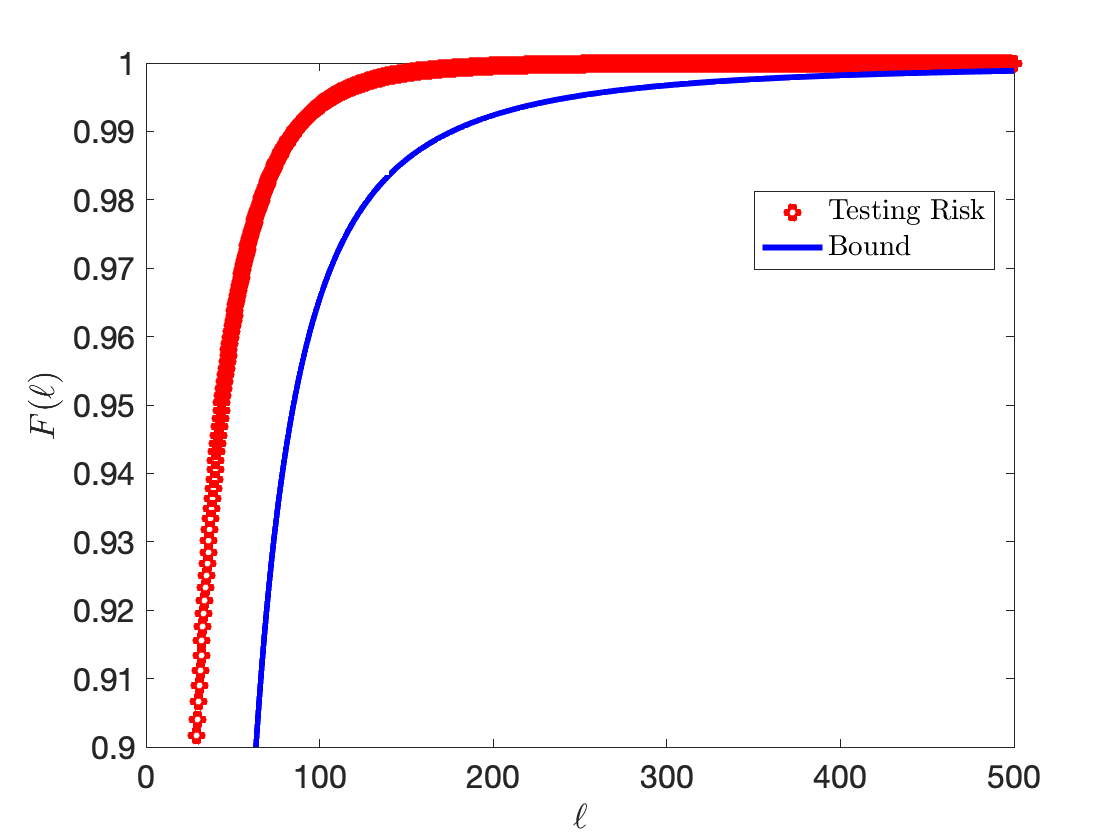}
	\includegraphics[width=0.33\textwidth]{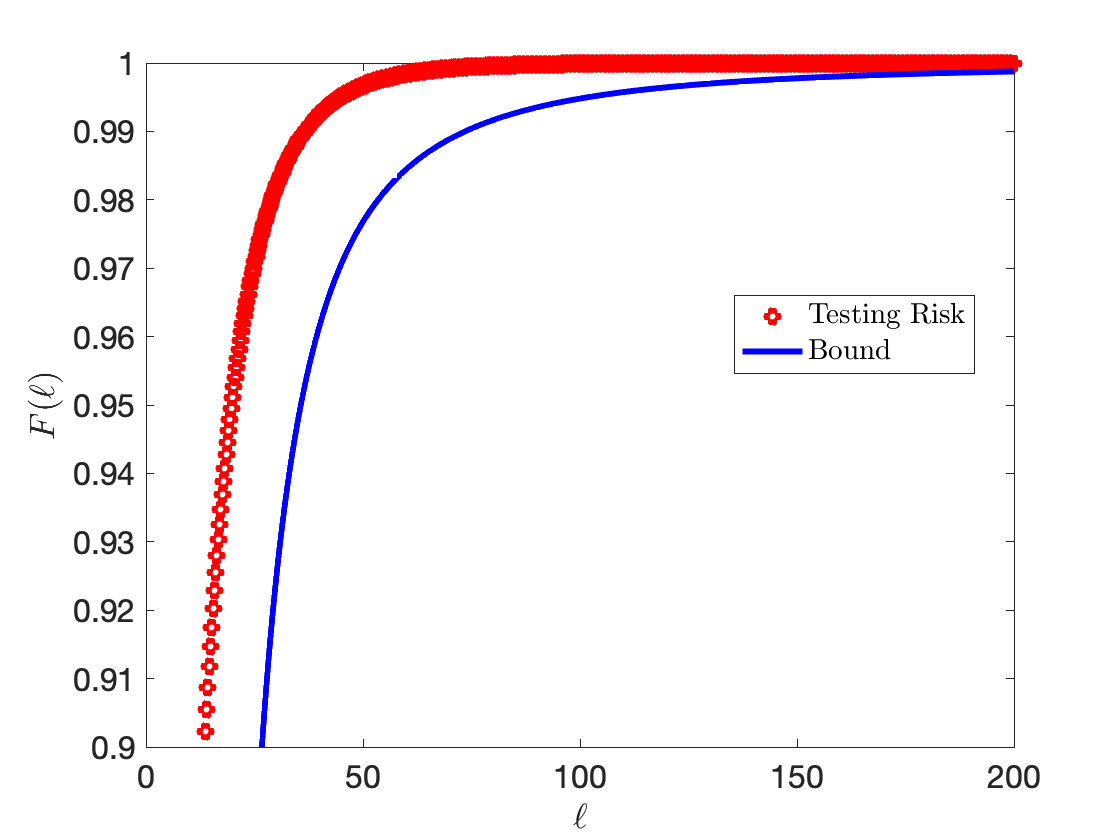}
		\includegraphics[width=0.33\textwidth]{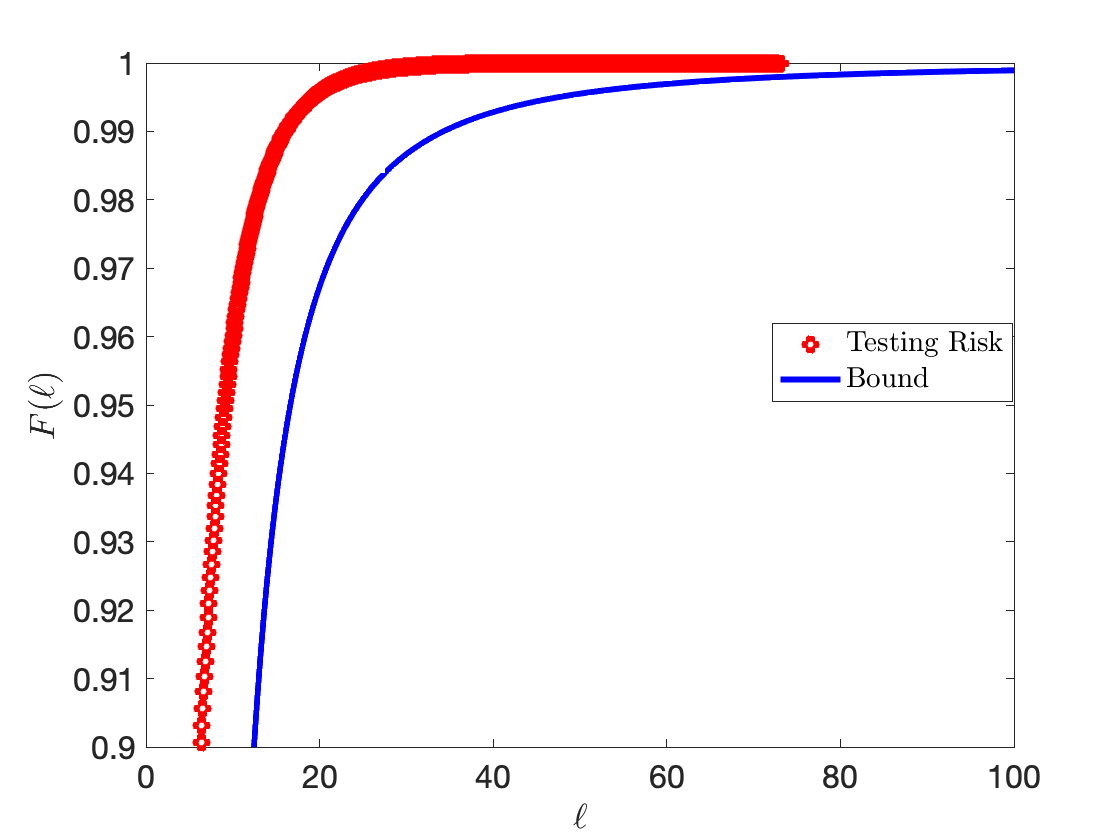}
	\caption{Empirical CDF and theoretical bounds for $n=40; m=1$. Left : N=8; Middle: N=20; Right: N=60}
	\label{fig:Risk2}
\end{figure}

Figure~\ref{fig:Risk2} verifies the empirical CDF (generated using $8 \times 10^6$ samples) for a case with signal to noise ratio of 1.

\begin{proof}
We will use the Chebyshev concentration inequality, which requires the variance. We first pursue the noise free case.  Following the setup of Proof 1, for the noiseless case
$$
\ell^2_{sys} \triangleq\hat{x}^T(I-a Q)A^T A(I-aQ)\hat{x} \hat{x}^T(I-a Q)A^T A(I-aQ)\hat{x}. $$
where

$a \triangleq \frac{\eta}{N}, Q \triangleq X X^T, A \triangleq W_1-W_0$ . 

\begin{align*}
\mathbb{E}_{\hat{x}}[\ell^2_{sys}] &= 2 \tr[A(I- a Q)(I- a Q)A^T A(I- a Q)(I- a Q)A^T] + (\tr[A(I- a Q)(I- a Q)A^T])^2\\
&= 2 \tr[B(I- a Q)(I- a Q)B(I- a Q)(I- a Q)] + (\tr[B(I- a Q)(I- a Q)])^2,
\end{align*}
where  $B \triangleq A^T A$.  We have used the following identity:

	If $X \sim \mathcal{N}(0,I)$ then
	$\mathbb{E}_X[x^T C^T C x x^T C^T C x] = 2\tr(CC^TCC^T)+ \tr[C C^T]^2.$ The manipulations below require expectations involving permutations of products of 8th order Gaussians. See Appendix.

\begin{align*}
\mathbb{E}[\ell^2_{sys}]&=  2\mathbb{E} \tr[B(I-2aQ+a^2QQ)B(I-2aQ+a^2QQ)]\\
&+ \mathbb{E}[\Tr[B(I-aQ)(I-aQ)]\Tr[B(I-aQ)(I-aQ)]]\\
&=  2\mathbb{E}\tr[BB-4aBBQ+2a^2BBQQ+4a^2BQBQ-4a^3BQBQQ\\
&+a^4BQQBQQ]+ \mathbb{E}[\tr[B]^2-4a\tr[B]\tr[BQ]+2a^2\tr[B]\tr[BQQ]+4a^2\tr[BQ]\tr[BQ]\\
&-4a^3\tr[BQ]\tr[BQQ]+a^4\tr[BQQ]\tr[BQQ]]\\
&=  2\tr[B^2](1-4aN+2a^2((n+1)N+N^2))+2\mathbb{E}\tr[4a^2BQBQ-4a^3BQBQQ\\
&+a^4BQQBQQ]+ \tr[B]^2(1-4aN+2a^2(1+n+N)) + \mathbb{E}[4a^2\tr[BQ]\tr[BQ\\
&-4a^3\tr[BQ]\tr[BQQ]+a^4\tr[BQQ]\tr[BQQ]]
\end{align*}
\begin{align*}
\mathbb{E}[\ell^2_{sys}]&=  2\tr[B^2](1-4aN+2a^2((n+1)N+N^2))\\
&+8 a^2 N((N+1)\tr[B^2]+\tr[B]^2)\\
&-8a^3 N ((4 + n + (3 + n) N + N^2) \tr[B^2] + (2 + n + 2 N)\tr[B]^2)\\
&+2 a^4 N((20 + n(11 + n) + 21 N + n (7 + n) N + 2 (3 + n) N^2 + N^3) \tr[B^2] \\
&+ (10 + 5 n + n^2 + 5 (2 + n) N + 4 N^2) \tr[B]^2)\\
&+\tr[B]^2(1-4aN+2a^2((n+1)N+N^2))\\
&+ 4a^2 N (N \tr[B]^2+2 \tr[B^2]) \\
&-4a^3 N(2(2 + n + 2 N) \tr[B^2] + (2 + N(1 + n + N))\tr[B]^2)\\
&+a^4 N(2(10 + 5n + n^2 + 5(2 + n)N + 4N^2)*\tr[B^2] + (n^2 N + 2 n (3 + N + N^2) \\
&+ (1 + N)(10 + N + N^2)) \tr[B]^2)\\
&=(1 + a N (-4 + 
a (2 (5 + n + 3 N) + 
a (2 + N) (3 + n + N) (-4 + a (5 + n + N)))))\\
& (2 \tr[B^2] + \tr[B]^2).
\end{align*}


The varjance $Var(\ell_{sys}) = \mathbb{E}[\ell^2_{sys}] - \mathbb{E}^2[\ell_{sys}]$ was verified using $10^{10}$ samples for several combinations of $N$ and $n$. Nevertheless the above formula is unwieldy. Examining the case where $m=1,w_0=0$ and $\eta=1$: $\mathbb{E}[\ell_{sys}] = \frac{n+1}{N}\norm{w_1-w_0}^2$ and 
$Var[\ell] = \left[ \frac{90+48 n+6 n^2}{N^3} + \frac{20+ 10 n +2 n^2}{N^2}\right]\norm{w_1-w_0}^4$.

 Chebyshev's inequality states that for any random variable $X$, and $\epsilon \in \reals^+$,

$$ P( |x-\mathbb{E}[X]| \geq \epsilon)  \leq  \frac{Var(X)}{\epsilon^2}. $$

Thus

\begin{align*}
P( |\ell_{sys}-\frac{n+1}{N}\norm{w_1-w_0}^2| \leq \epsilon)   &\geq 1- \frac{1}{\epsilon^2} \left[ \frac{90+48 n+6 n^2}{N^3} + \frac{20+ 10 n +2 n^2}{N^2}\right]\norm{w_1-w_0}^4 \\
\end{align*}
Define $\delta \triangleq \frac{Var[\ell]}{\epsilon^2} $. Therefore with at least a probability of $1-\delta$.

\begin{equation}
\ell_{sys} \leq \left[\frac{n+1}{N}+ \sqrt{\frac{ \left[ \frac{90+48 n+6 n^2}{N^3} + \frac{20+ 10 n +2 n^2}{N^2}\right]}{\delta}} \right] \norm{w_1-w_0}^2.
\end{equation}

Now we consider the noise terms.
Define $D \triangleq A(I-aQ)  - a \sigma Z X^T$
\begin{align*}
\ell^2 &=  (D \hat{x} + \sigma \hat{z})^T (D \hat{x} + \sigma \hat{z}) (D \hat{x} + \sigma \hat{z})^T (D \hat{x} + \sigma \hat{z}) \\
\mathbb{E}_{\hat{x}}[\ell^2] &=2 \tr[D D^T D D^T] + 4 \sigma^2 \hat{z}^T D D^T \hat{z} + (\tr[D D^T]+\sigma^2 \hat{z}^T \hat{z})^2 \\
&=2 \tr[D D^TD D^T] + 4 \sigma^2 \hat{z}^T D D^T \hat{z} + \tr[D D^T]^2+\sigma^4 \hat{z}^T \hat{z} \hat{z}^T \hat{z} + 2 \sigma^2 \hat{z}^T \hat{z} \tr[D D^T]\\
\mathbb{E}_{\hat{z} \hat{x}}[\ell^2] &=2 \tr[D D^T D D^T] + 4 \sigma^2  \tr[D D^T] + \tr[D D^T]^2+(2m+m^2) \sigma^4 + 2 \sigma^2 m \tr[D D^T].
\end{align*}

Let's switch to $m=1$ for tractability
$  c^T  \triangleq A(I-aQ) ; C^T \triangleq  - a \sigma X^T$ and $z = Z^T$. The expressions below require new identities that are presented in the Appendix.

\begin{align*}
\mathbb{E}_{\hat{z} \hat{x}}[\ell^2] &= 3 (Cz+c)^T(Cz+c)(Cz+c)^T(Cz+c)+ 6 \sigma^2  (Cz+c)^T(Cz+c) + 3 \sigma^4 \\
\mathbb{E}_{z \hat{z} \hat{x}}[\ell^2] &= 6 \Tr[C C^T C C^T] + 12 c^T C C^T c + 3 (\Tr[CC^T]  + c^T c)^2 + 6 \sigma^2 (\tr[C C^T]+c^T c) + 3 \sigma^4 \\
&= 6 a^4 \sigma^4\Tr[Q Q] + 12a^2 \sigma^2 \Tr[BQ-2aBQQ +a^2BQQQ] +3 a^4 \sigma^4 \Tr[Q]\Tr[Q]{\color{blue}+ \mathbb{E}[\ell_{sys}^2]}\\
&+6a^2\sigma^2 \Tr[Q] (A A^T -2aAQA^T+a^2 AQQA^T) + 6 \sigma^4 a^2 \Tr[Q] \\
&+ 6 \sigma^2 \tr[B-2aBQ+a^2 BQQ] + 3 \sigma^4\\
\mathbb{E}[\ell^2]&= 6 a^4 \sigma^4 Nn (N+n+1) + 12a^2 \sigma^2 \Tr[B](N-2aN(N+n+1)\\
& +a^2N (4 + n^2 + 3 n (1 + N) + N (3 + N))) +3 a^4 \sigma^4 N (2n +N n^2)\\    
&+6a^2\sigma^2 \tr[B] (Nn  -2aN (2+Nn)+a^2 N(1+n+N)(4+nN)) \\
&+ 6 \sigma^4 a^2 nN + 6 \sigma^2 \tr[B](1-2aN+a^2 N(N+n+1)) + 3 \sigma^4  {\color{blue}+ \mathbb{E}[\ell_{sys}^2]}
\end{align*}

\noindent For $\eta = 1$, we get
\begin{align*}\mathbb{E}[\ell^2]&= \sigma^4 \left( 3 + \frac{12 n+ 6 n^2 }{N^3} + \frac{6 n + 3 n^2}{N^2} + \frac{6 n}{N}\right) \\ 
& + \tr[B]  \sigma^2\left(  \frac{72   + 60 n  +12 n^2 }{N^3} + \frac{12  + 18 n +6 n^2 }{N^2} + \frac{6 n}{N} + \frac{6}{N} \right)  {\color{blue}+ \mathbb{E}[\ell_{sys}^2]}.
\end{align*}

\begin{align*}
Var[\ell] &=\sigma^4\left(2  + \frac{12 n+6 n^2}{N^3} + \frac{6 n + 2 n^2}{N^2} + \frac{4 n}{N}\right)\\
&+\sigma^2 \tr[B]\left(\frac{72  +  60 n  + 12 n^2}{N^3} + \frac{12 + 16 n + 4 n^2}{N^2} + \frac{4 + 4 n}{N}\right) \\
&+\tr[B]^2 \left( \frac{90  + 48 n + 6 n^2}{N^3} + \frac{20  + 10 n + 2 n^2}{N^2} \right).
\end{align*}
\end{proof}

\section{Results for Least Norm and Least Squares Regression}
\label{sec:LSQ}

\begin{theorem}[Properties of Least Squares Regression]\label{thm:LSQ}
Given $X \triangleq [x_1 \ \ x_2  \ \ ...\ \ x_N] \in \reals^{n \times N}$, which is a i.i.d standard normal random matrix and $Y \triangleq [y_1 \ \ y_2 \ \ .. \ \ y_N] \in \reals^{1 \times N}$ with $y_i = w_1^T x_i + \sigma^2 z_i$ and $z_i \in \reals^m$ is an i.i.d. sample from the standard normal distribution. Then, the generalization error in Least norm (centered on $W_0$) and Least Squares regression has the following properties

\begin{align*}
\mathbb{E}[\ell] &=\norm{W_1 - W_0}^2 \left(1-\frac{N}{n} \right) { + \sigma^2 \left(1 + \frac{N}{n-N-1} \right)} \ \  ; \ \ 1 \leq N \leq n-1 \\
\mathbb{E}[\ell] &= \sigma^2 \left(1 + \frac{n}{N-n-1} \right) \ \  ; \ \ n+1 \leq N  \\
\mathbb{E}[\ell^2] &= \sigma^4 \frac{3 (N-1) (N-3)}{(N-n-1) (N-n-3)} \ \ ; N  > n+3 \\
\mathbb{E}[\ell^2] &=\sigma^4 \frac{3 (n-1) (n-3)}{(n-N-1) (n-N-3)} + 6 \sigma^2 \norm{W_1 - W_0}^2 \frac{(n-1)(n-N)}{n(n-N-1)} \\
& + 3 \norm{W_1 - W_0}^4  \frac{(n-N)(n-N+2)}{n(n+2)}  \ \ ; \ \ n \ge N+3.
\end{align*}
\end{theorem}

\begin{figure}
	\centering
	\includegraphics[width=0.75\textwidth]{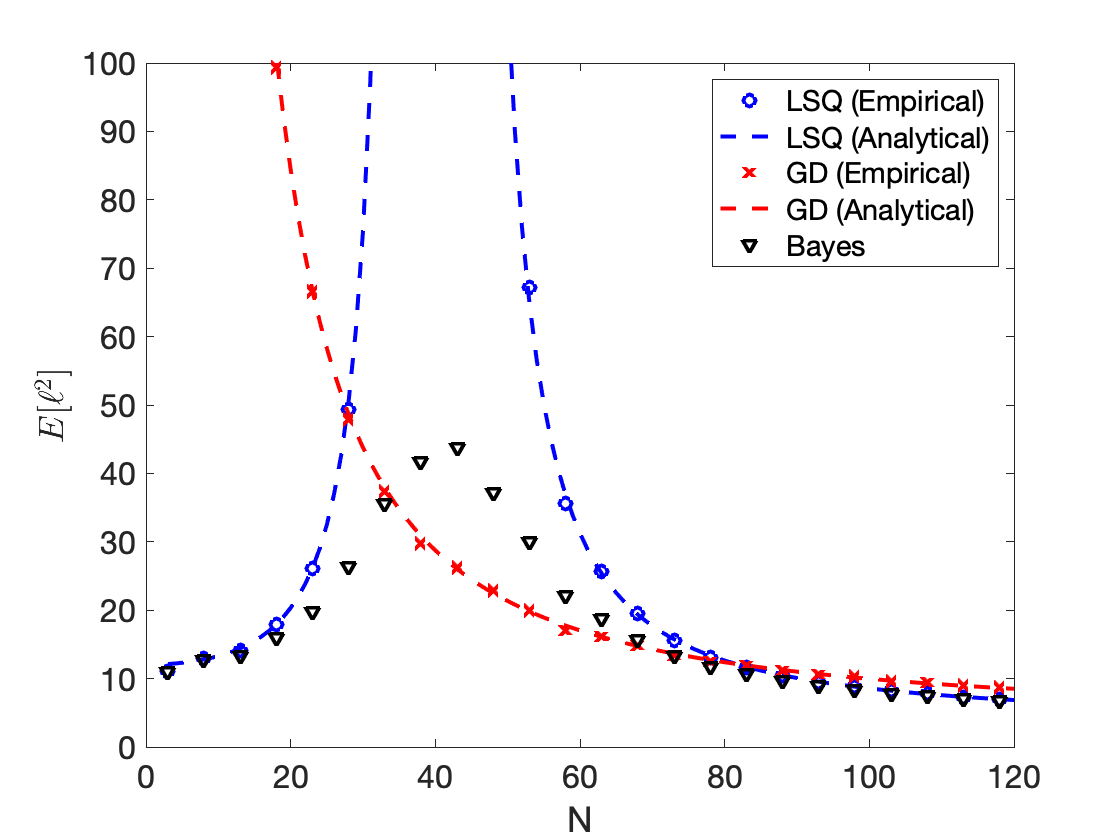}
	\includegraphics[width=0.49\textwidth]{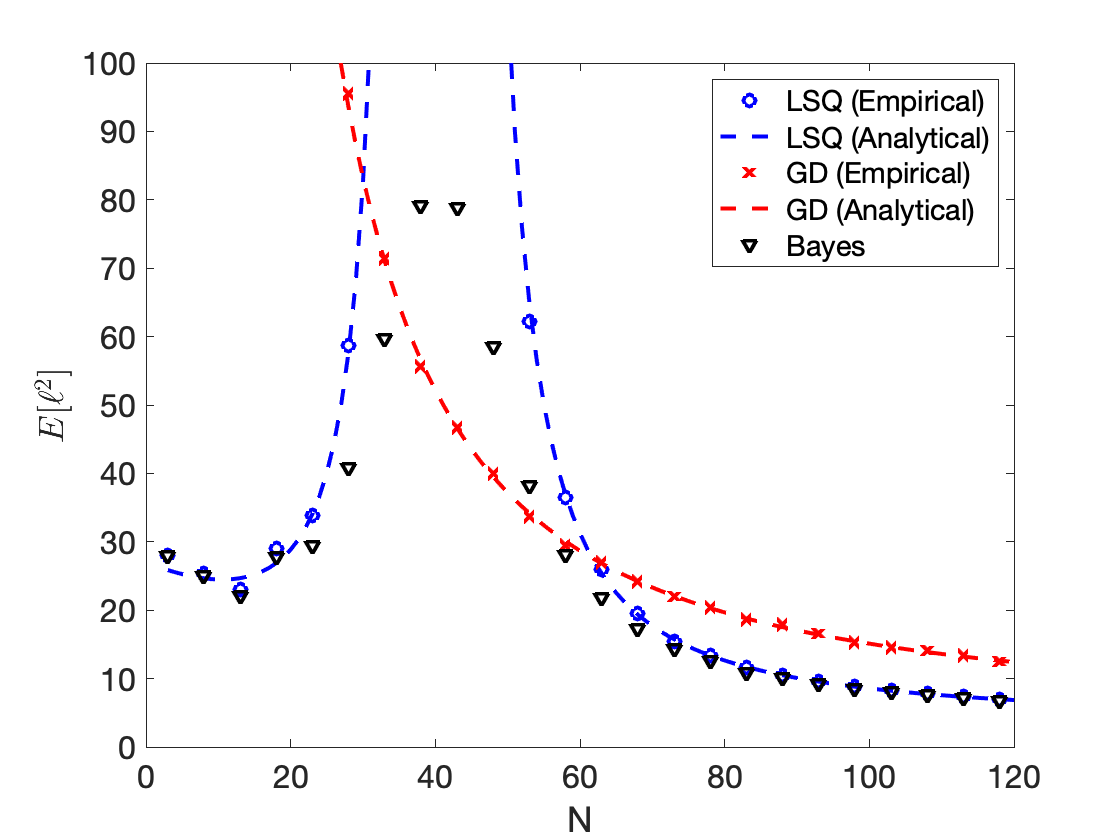}
	\includegraphics[width=0.49\textwidth]{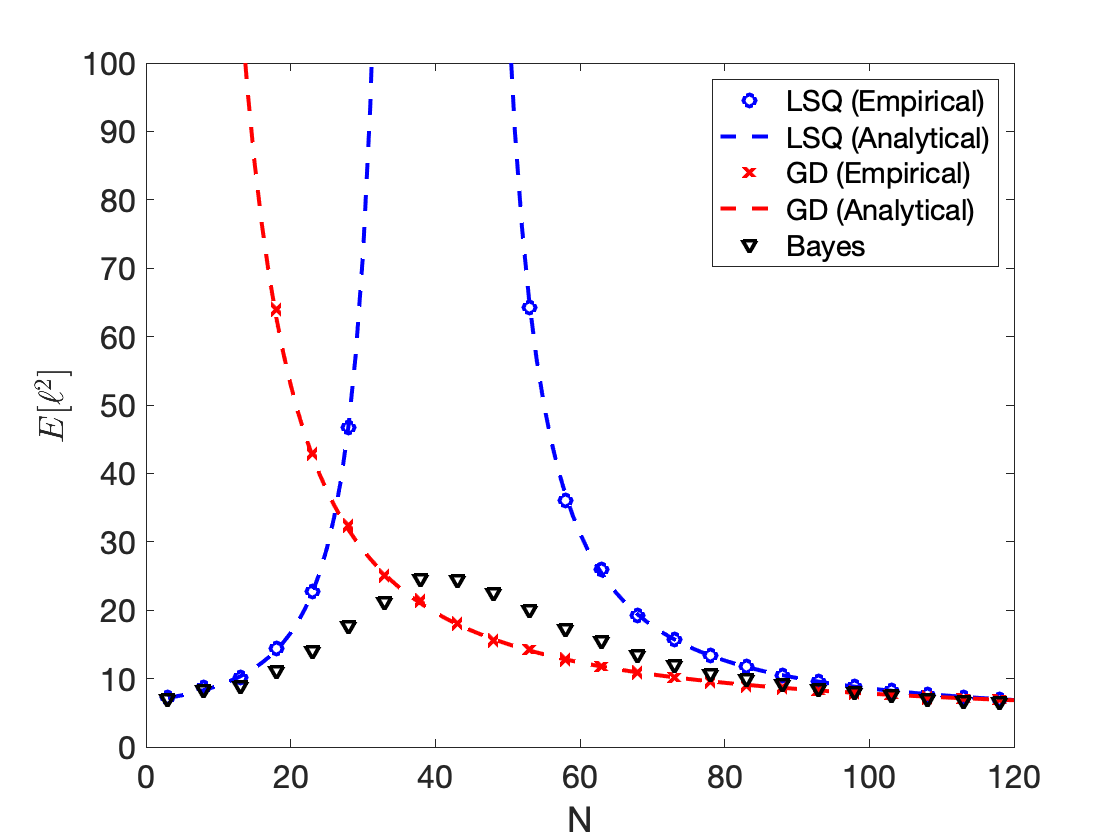}
	\caption{$\mathbb{E}[\ell^2]$ for $n=40; m=1; \eta=1$. Top : $\norm{W_1}^2/\sigma^2=1$; Bottom left : $\norm{W_1}^2/\sigma^2=2$; Bottom right: $\norm{W_1}^2/\sigma^2=0.5$ }
	\label{fig:EL2}
\end{figure}

\begin{figure}
	\centering
	\includegraphics[width=0.33\textwidth]{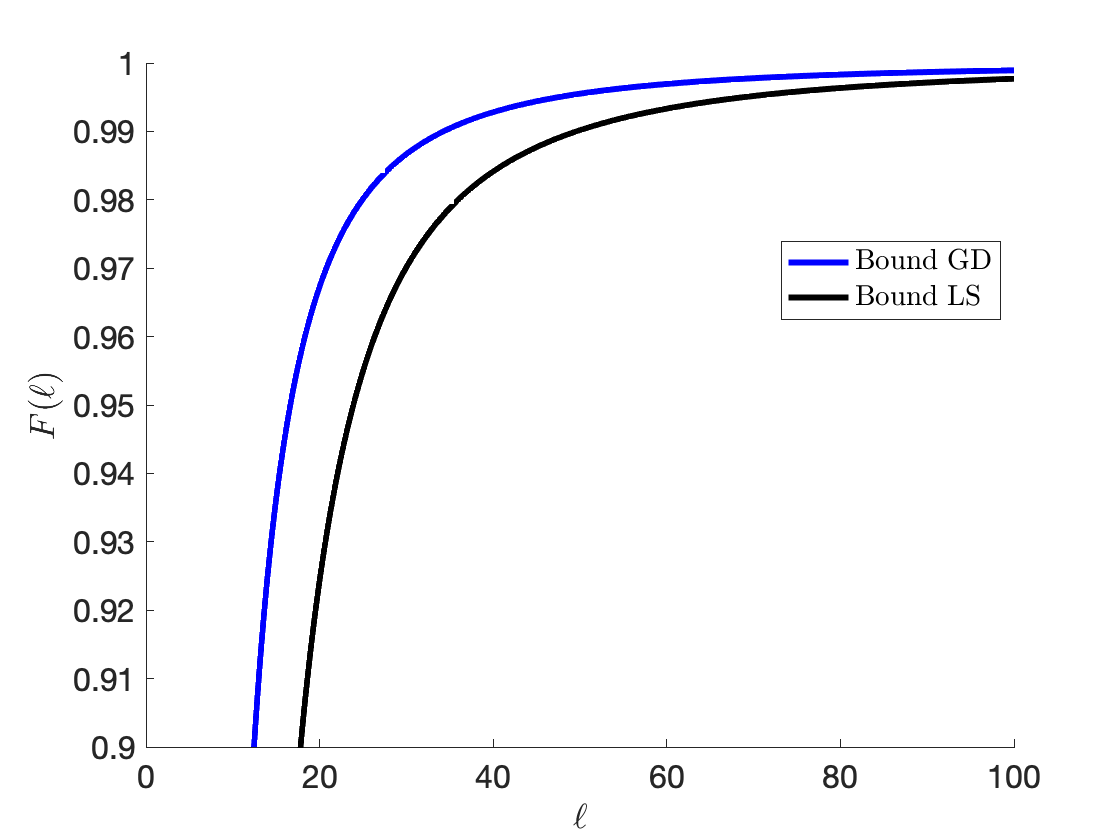}
	\includegraphics[width=0.33\textwidth]{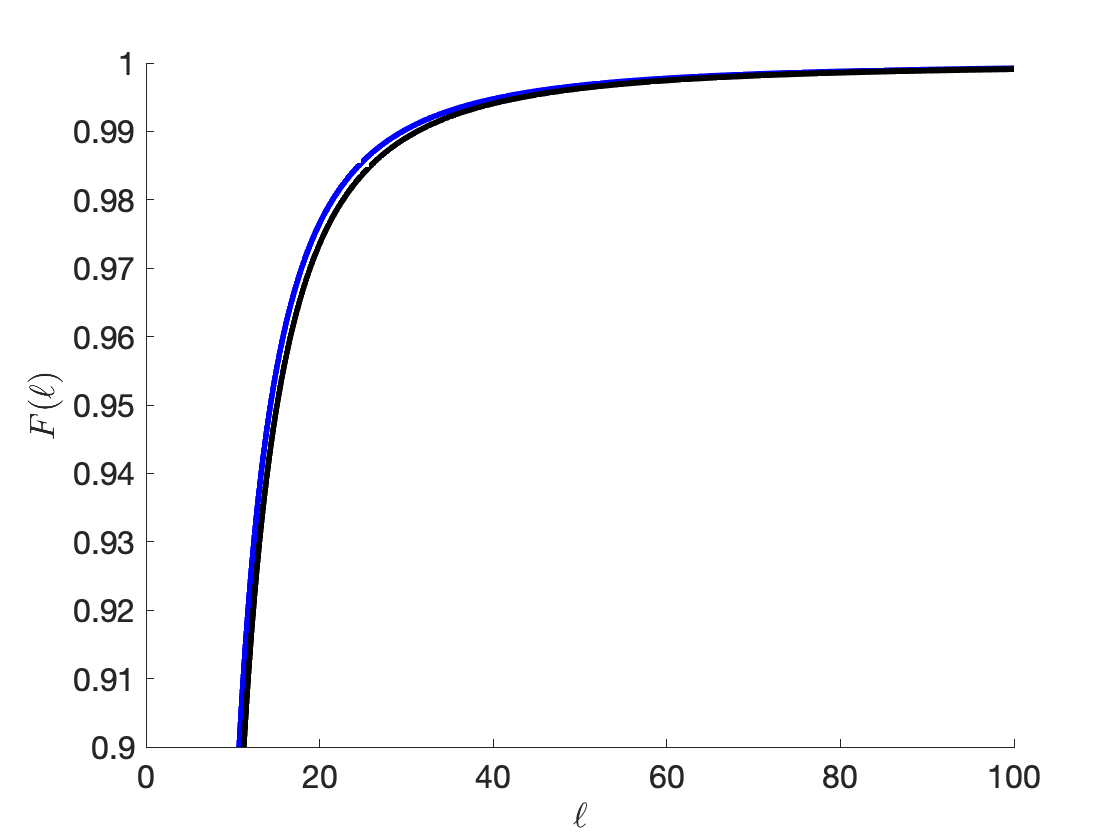}
	\caption{Analytical bounds for $n=40; m=1$. Left : N=60;  Right: N=80}
	\label{fig:Risk3}
\end{figure}

\noindent {\em Note:} A version of the first two equalities can be found in Belkin et al.~\cite{belkin2020two} and are verified in Figure~\ref{fig:EL}. The last two equalities above are novel to the knowledge of the author, and are verified in Figure~\ref{fig:EL2}. Hastie et al.~\cite{hastie2022surprises} have derived formulae  for variance under slightly more general conditions, but these results are asymptotic in nature. Very recently, Zhou et al.~\cite{zhou2024optimistic} derive a related expression, but they do not assume $\hat{x}$ and $\hat{y}$ are random as in the present context~\footnote{They have cited an  unpublished, early version of this work~\cite{duraisamy2021basic}}.  Belkin et al.~\cite{belkin2020two} derive statistical bounds that are  related. 
These results can be employed with concentration inequalities as in the previous section to construct bounds for the generalization error.  Figure ~\ref{fig:Risk3} shows such a comparison for the $n=40$ example considered in previously.

\begin{proof}
	
\noindent {\bf \underline{Over-parametrized case}}

	Let us consider the over-parametrized case first:
	\begin{align*}
	\ell &= ||(W_1 - W)\hat{x} {+ \sigma \hat{z} }||_2^2 = ||  \sigma \hat{z}  - \sigma Z X^+ \hat{x}  ||_2^2\\
	\mathbb{E}[\ell]&= \mathbb{E}[(\sigma \hat{z}  - \sigma Z X^+ \hat{x} )^T(\sigma \hat{z}  - \sigma Z X^+ \hat{x} )] \\&=  \mathbb{E}[\sigma^2 \hat{z}^T \hat{z}]   +\mathbb{E}[\sigma^2 \hat{x}^T X^{+T} Z^T Z  X^+ \hat{x}] \\
	&=  \sigma^2    + \sigma^2 \mathbb{E}_{X Z}[ \tr[X^{+T} Z^T Z  X^+]] \\
	&=  \sigma^2    + \sigma^2 \mathbb{E}_{X}[ \tr[X^{+T}   X^+]] \\
	&=  \sigma^2    + \sigma^2 \mathbb{E}_{X}[ \tr[(X X^T)^{-1}]] \\
	&= \sigma^2 \left(1 + \frac{n}{N-n-1} \right).
	\end{align*}

	\begin{align*}
	\ell^2 &= (\sigma \hat{z}  - \sigma Z X^+ \hat{x} )^T(\sigma \hat{z}  - \sigma Z X^+ \hat{x} )(\sigma \hat{z}  - \sigma Z X^+ \hat{x} )^T(\sigma \hat{z}  - \sigma Z X^+ \hat{x} )\\
	\mathbb{E}_{\hat{x}}[\ell^2] &= \sigma^4 ( 2 \tr[Z X^+ X^{+T} Z^T Z X^+ X^{+T} Z^T] + 4  \hat{z}^T Z X^+ X^{+T} Z^T \hat{z} + \tr[Z X^+ X^{+T} Z^T ]^2 + \hat{z}^T \hat{z} \hat{z}^T \hat{z}\\
	& + 2  \hat{z}^T \hat{z} \tr[Z X^+ X^{+T} Z^T ])\\
	\mathbb{E}_{\hat{z}\hat{x}}[\ell^2]&=  \sigma^4 ( 2 \tr[z^T (X X^T)^{-1} z z^T (X X^T)^{-1} z] + 4  \tr[ z^T (X X^T)^{-1} z] + \tr[z^T (X X^T)^{-1} z ]^2 + 3 \\
	&+ 2  \tr[z^T (X X^T)^{-1} z ])\\
	\mathbb{E}_{z\hat{z}\hat{x}}[\ell^2]&=  \sigma^4 ( 6 \tr[(X X^T)^{-2}] + 3 \tr[X X^T)^{-1}]^2 + 6  \tr[(X X^T)^{-1}] + 3)\\
	&\sim \sigma^4( 6 \tr(W_n^{-1}(I,N)W_n^{-1}(I,N)) + 3 \tr(W_n^{-1}(I,N))\tr(W_n^{-1}(I,N))+ 6 \tr( W_n^{-1}(I,N))  +3),
	\end{align*}
	where $W_m^{-1}$ represents the Inverse Wishart distribution. The expectation of the last expression can be compactly reduced using standard Wishart distribution identities, except the term involving the trace of the product of inverse Wishart matrices, which we obtain from Pielaszkiewicz \& Holgersson~\cite{pielaszkiewicz2019mixtures} (page 8). With this identity,
	\begin{align*}
	\mathbb{E}_X \mathbb{E}_{\xh} \mathbb{E}_{z} \mathbb{E}_{\hat z} \ell^2&= 6 \sigma^4 \frac{(N-1)n}{(N-n-3)(N-n-1)(N-n)} +3 \sigma^4 \frac{n(n(N-n-2)+2)}{(N-n-3)(N-n-1)(N-n)}\\
	& + 6 \sigma^4 \frac{n}{N-n-1}  + 3 \sigma^4 \\
	&=\sigma^4 \frac{3 (N-1) (N-3)}{(N-n-1) (N-n-3)},
	\end{align*}
	as long as $N>n+3$.

	\noindent {\bf \underline{Under-parametrized case}}
	
Next, we consider the under-parameterized case for which the Lagrangian \( \mathcal{L} \) is defined as:
	\[
	\mathcal{L}(W, \Lambda) = \text{tr}((W - W_0)(W - W_0)^T) + \text{tr}(\Lambda^T (W X - Y))
	\]
	where \( \Lambda \) is a  matrix of Lagrange multipliers.	
	It is easily shown that first order optimality yields $W=W_0 + (Y - W_0 X) (X^T X)^{-1} X^T$. In the present setup, $Y = W_1 X {+ \sigma Z}$. Therefore
$
W = W_0(I-P) + W_1 P { + \sigma Z X^+} ,
$
where $X^+ \triangleq (X^T X)^{-1} X^T$ and  $ P  \triangleq  X (X^T X)^{-1} X^T$.

For a test location $\hat{x}$,
$$ \ell =  ||(W_1 - W)\hat{x} {+ \sigma \hat{z}}||_2^2 = || (W_1-W_0)(I-P)  \hat{x} {+ \sigma \hat{z} - \sigma Z X^+ \hat{x}}||_2^2$$

Taking the expectation of the first term (wrt) $\hat{x}$, we have  $  \tr[(I-P)(W_1-W_0)^T(W_1-W_0)(I-P)] = \tr[(I-P)(W_1-W_0)^T(W_1-W_0)]$
Then taking the expectation of the first term wrt to X, we have
\begin{align*}
 \tr[(W_1-W_0)^T(W_1-W_0)\mathbb{E}[I-P]] = \tr[(W_1-W_0)^T(W_1-W_0)]\left(1-\frac{N}{n}\right).
\end{align*}
The contribution of the second term to $\mathbb{E}[\ell]$ is $\sigma^2$.

The third term is $
\mathbb{E}_{\hat x}[\hat{x}^T X^{+T} Z^T Z X^+ \hat{x}] = \Tr[X^{+T} Z^T Z X^+]$, thus
$\mathbb{E}_Z[\Tr[X^{+T} Z^T Z X^+]] = \Tr[X^{+T} X^+] = \Tr[(X^T X)^{-1}]$. Finally,
$$\mathbb{E}_X[\tr[(X^T X)^{-1}]] = \frac{N}{n-N-1}.$$

Therefore $$\mathbb{E}[\ell] = ||W_1 - W_0||^2 \left(1-\frac{N}{n} \right) { + \sigma^2 \left(1 + \frac{N}{n-N-1} \right)}. $$

Now we consider the second moment. Define $C  \triangleq  -  \sigma  X^{+T}; c\triangleq (I-P) b; z \triangleq Z^T; b \triangleq (W_1-W_0)^T.$
\begin{align*}
\mathbb{E}_{z \hat{z} \hat{x}}[\ell^2] &= 6 \Tr[C C^T C C^T] + 12 c^T C C^T c + 3 (\Tr[CC^T]  + c^T c)^2 + 6 \sigma^2 (\tr[C C^T]+c^T c) + 3 \sigma^4 \\
&= 6 \sigma^4 \Tr[X^{+T} X^{+} X^{+T} X^{+}] + 12 \sigma^2 b^T (I-P)  X^{+T} X^{+} (I-P)b + 3 \sigma^4 \Tr[X^{+T}X^{+}]^2  \\
&+ 3 b^T (I-P)b b^T (I-P)b  \\
&+ 6 \sigma^2 \tr[ X^{+T}  X^{+}]b^T(I-P)b + 6 \sigma^4 \tr[X^{+T} X^{+}]+6 \sigma^2 b^T(I-P)b+ 3 \sigma^4 \\
&= 6 \sigma^4 \Tr[(X^T X)^{-2}] + 3 \sigma^4 \Tr[(X^T X)^{-1}]^2  + 3 b^T (I-P)b b^T (I-P)b  \\
&+ 6 \sigma^2 \tr[ (X^T X)^{-1}]b^T(I-P)b + 6 \sigma^4 \tr[(X^T X)^{-1}]+6 \sigma^2 b^T(I-P)b+ 3 \sigma^4 \\
\mathbb{E}[\ell^2] &=6 \sigma^4 \frac{(n-1)N}{(n-N-3)(n-N-1)(n-N)} \\&+ 3 \sigma^4 \frac{N(N(n-N-2)+2)}{(n-N-3)(n-N-1)(n-N)}  + 3 \norm{b}^4  \frac{(n-N)(n-N+2)}{n(n+2)}  \\
&+ 6 \sigma^2 \frac{N}{n-N-1}\norm{b}^2 \left(1-\frac{N}{n}\right) + 6 \sigma^4\frac{N}{n-N-1}+6 \sigma^2 \norm{b}^2 \left(1-\frac{N}{n}\right)+ 3 \sigma^4\\
&=\sigma^4 \frac{3 (n-1) (n-3)}{(n-N-1) (n-N-3)} + 6 \sigma^2 \frac{(n-1)(n-N)}{n(n-N-1)}\norm{b}^2  + 3 \norm{b}^4  \frac{(n-N)(n-N+2)}{n(n+2)}. 
\end{align*}

In the derivation above, we have used the following 
\begin{enumerate}
\item Noting that $I-P$ will have $n-N$ unity eigenvalues (rest are zero). We will write $(I-P) = \sum_{i=1}^{n-N} v_i v_i^T$, where $v_i$ are its eigenvectors.
\begin{align*}
b^T(I- P)b b^T(I- P)b&= b^T \sum_{i=1}^{n-N} (v_i v_i^T) b b^T \sum_{j=1}^{n-N} (v_j v_j^T) b\\
&=(n-N)b^T  v_1 v_1^T b b^T  v_1 v_1^T b +  (n-N)(n-N-1)b^T  v_1 v_1^T b b^T  v_2 v_2^T b\\
&= (n-N) (b^T  v_1)^4  +  (n-N)(n-N-1) (b^T  v_1)^2 (b^T  v_2)^2\\
\mathbb{E}[b^T(I- P)b b^T(I- P)b]&= \norm{b}^4 (n-N) \left( \frac{3}{n(n+2)}  +  (n-N-1) \frac{1}{n(n+2)}  \right)\\
&= \norm{b}^4  \frac{(n-N)(n-N+2)}{n(n+2)}.
\end{align*}
The penultimate step uses identities 29 and 30 in Section~\ref{sec:identities}.

\item $X^+ P = X^+$ and hence $X^+(I-P) = 0.$

\item Write $X = U \Sigma V^T$
 \begin{align*} 
\tr[ (X^T X)^{-1}]b^T(I-P)b & = \left(\tr[(V \Sigma U^T U \Sigma V^T )^{-1}]  \right) \left(b^T (I-U\Sigma V^T (V \Sigma U^T U \Sigma V^T)^{-1} V \Sigma U^T)b\right)\\
&= \left(\sum_{i=1}^{N} \frac{1} {\sigma_i^2} \right) \left( b^T(I-UU^T)b\right).
\end{align*}
Since the elements of $X$ are i.i.d Gaussian random variables, the singular vectors $U$ and $V$ 
are uniformly distributed on the unit sphere. Because of the rotational invariance of the singular vectors, the columns of 
$U$ and $V$  do not depend on the particular values of the singular values $\sigma$, and thus the two bracketed terms above are independent of each other. Thus
 \begin{align*} 
\mathbb{E}[\tr[ (X^T X)^{-1}]b^T(I-P)b] & = \mathbb{E}[\tr[ (X^T X)^{-1}]] \mathbb{E}[b^T(I-P)b] \\
&= \frac{N}{n-N-1}  \norm{b}^2 \left(1-\frac{N}{n} \right).
\end{align*}

\end{enumerate}

\end{proof}

\begin{figure}[!htb]
	\centering
	\includegraphics[width=0.75\textwidth]{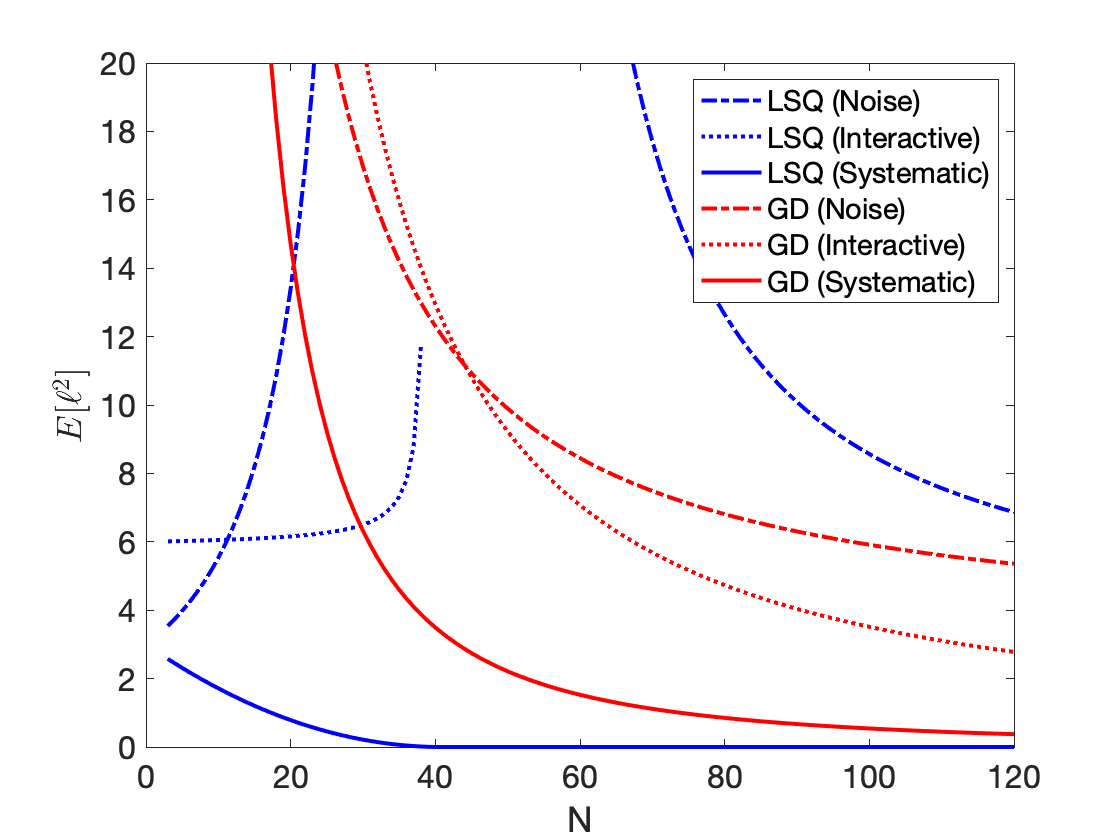}
	\caption{$\mathbb{E}[\ell^2]$ for $n=40; m=1$ . Break down shows pure noise, interactive and systematic components for Least Squares and one-shot gradient descent (with $\eta=1$). }
	\label{fig:EL2b}
\end{figure}

Figure~\ref{fig:EL2b} shows the different components of $\mathbb{E}[\ell^2]$ for a sample case, comparing the systematic (terms involving $\norm{W_1-W_0}^4$), pure noise (terms involving $\sigma^4$) and 'interactions' (terms involving $\norm{W_1-W_0}^2 \sigma^2$). Even for the un-optimized gradient descent (i.e. $\eta=1$), the implicit regularization appears to be highly beneficial.

\section{Connections to existing  work \& Conclusions}
\label{sec:conclusions}
This study investigated the performance of gradient descent in a in-context linear regression setting, providing qualitative insights and quantitative characterization of the statistical properties of the generalization error. The derived generalization error bounds were contrasted with those from least squares regression, demonstrating that a single step of gradient descent can offer a comparable performance in certain contexts, especially in high noise settings.  An expression for the optimal step size was derived. The  analysis of systematic and noise components offered a comprehensive view of the factors contributing to generalization error. 

The fact that gradient descent can provide reasonable results with a single step has implications for reducing computational complexity, especially in one-shot scenarios and resource-constrained environments. Additionally, generalization error bounds were extended to least squares and least norm regression. This work uncovered new identities involving high-order products of Gaussian random matrices, which may have broader applications in regression tasks and beyond. All of the derived results are verified using empirical computations on a large number of samples. It is  intriguing that in much of the literature, probabilistic bounds are often not verified using numerical experiments (indeed, there are clear exceptions, for instance~\cite{hastie2022surprises,mei2019generalization,zhou2024optimistic}). 

Our study, while restricted in scope to well-specified models, addresses some key gaps in the literature: Despite extensive research on convergence and optimality in gradient-based methods, finite sample, non-asymptotic results that do not rely on arbitrary constants are rare in the literature, and in the case of probabilistic bounds, this appears to be the case for linear regression as well. Closely related to our work -- though focused more on Transformers directly -- Mahankali et al.~\cite{mahankali2023one} show that  linear self-attention layers can emulate one step of gradient descent on least-squares regression tasks, and exhibit optimal in-context learning capabilities in synthetic linear regression scenarios. Zhang et al~\cite{zhang2023trained} show that transformers emulate gradient descent by achieving global minimums through suitable initialization, enabling competitive prediction errors on new tasks, while being sensitive to covariate shifts. These works are focused on optimality and do not consider bounds.


Existing literature on gradient descent and both linear and non-linear regression is extensive and rigorous, and might prompt questions about the value of a simpler setting of the  present work.  Nevertheless,  clean, finite sample results without arbitrary constants can yield much insight into the behavior of more complex architectures and regression tasks. As a comparable example,  Belkin et al.~\cite{belkin2020two} also consider well-specified linear regression and explain the so-called double-descent phenomenon which challenges classical notions of the bias-variance trade-off. This finding has been reproduced in more complex problems  involving deep neural networks (e.g. ~\cite{nakkiran2021deep}), and thus has led to profound implications for the design and understanding of learning algorithms.

Overall, this study underscores the potential for a single step of gradient descent to generalize effectively in in-context learning scenarios. Future research could extend these findings to more complex regression tasks, including non-linear and incomplete parametrizations, i.e. considering model form errors. Additionally, the implications of these results can be explored on transformer-based architectures and other machine learning algorithms on practical applications.

\section{Appendix}

\subsection{Identities}
\label{sec:identities}
Some of the derivations presented in the manuscript were extremely lengthy. To aid further work, we present a compact set of identities below.  A few of these can be easily derived using the excellent Matrix reference manual~\cite{matrixref}. Many of the below expressions, however, require many hours of manipulations. Proofs are also provided for a few of the more complex identities.

\noindent All the expressions assume that $x_i \sim \mathcal{N}(0,I_n)$  ;  $Q\triangleq \sum_{i=1}^N x_i x_i^T$ ; $B \in \reals^{n \times n}$ ; $b \in \reals^n$  

\begin{enumerate}
	
		
	\item  $ \mathbb{E}[x x^T  x x^T ] = (2+n) I$
	
	
	\item $\mathbb{E}[\tr[BQ]] = N \tr[B]$
	
		\item  $ \mathbb{E}[x x^T  x x^T x x^T] = (8+6n+n^2) I$
	
	\item  $ \mathbb{E}[x x^T B  x x^T ] =  B + B^T + \tr[B] I$
	
	\item  $ \mathbb{E}[x^T B  x x^T B x] =  \tr[B (B+B^T)] + \tr[B]^2$

		\item $\mathbb{E}[\tr[Q]\tr[Q]] = N(2n+N n^2)$
		
		\item $\mathbb{E}[\tr[QQ]] = N n (N+n+1)$
		
	The following 3 identities~\cite{pielaszkiewicz2019mixtures} assume invertibility of $X X^T$ which is almost surely guaranteed under the present settings for $N>n$. Switch n and N for analogous identities involving the inverse of $X^T X$. 	
	\item  $\tr[(XX^T)^{-1}] = \frac{n}{N-n-1}$

	\item $\tr[(XX^T)^{-2}] = \frac{(N-1)n}{(N-n-3)(N-n-1)(N-n)}$
	
	\item $\tr[(XX^T)^{-1}]^2 = \frac{n(n(N-n-2)+2)}{(N-n-3)(N-n-1)(N-n)}$

			\item $\mathbb{E}[x^T x b^T x x^T x x^T b] = b^T b (n^2+6n+8)$

		\item 
			$\mathbb{E}[(Bx + b)^T(Bx + b) (Bx + b)^T(Bx + b)] = 2\tr(BB^TBB^T) + 4b^TBB^Tb + (tr(BB^T) + b^T b)^2.$

			\item $\mathbb{E}[\tr[Q]b^TQQb]] = b^T b N (1 + n + N) (4 + n N)$
			
						\item $\mathbb{E}[\tr[Q]b^TQb]] = b^T b N (2 + n N) $
			
			\item $\mathbb{E}[QQQ] = N (4 + n^2 + 3 n (1 + N) + N (3 + N)) I $

			\item $\mathbb{E}[\tr[BQQ]] = \tr[B] N(N+n+1)$
		\item $\mathbb{E}[\tr[Q]b^TQb] = N(2+n N)b^T b$
				\item $\mathbb{E}[\tr[Q]b^TQ Q b] = N (1+n+N) (4+nN)b^Tb$
				
					\item $ \mathbb{E}[x^T B  x b^T x x^T] =  b^T(B+B^T+\tr[B]I)$

	The following assume that B is symmetric
		
	\item  $ \mathbb{E}[x^T B  x x^T B x x^T x] =  (n+4) \tr[B]^2+ (8+2n) \tr[B^2]$ 
	
		\item  $ \mathbb{E}[x^T B  x x^T x x^T B x x^T x] =  (n^2+24+10 n) \tr[B]^2+ (2n^2+20n+48) \tr[B^2]$
	
	\item $\mathbb{E}[\tr[BQBQ]] = N ((N+1)\tr[B^2]+\tr[B]^2)$
	
		\item $\mathbb{E}[\tr[BQBQQ]] = N((4 + n + (3 + n)N + N^2)\tr[B^2] + (2 + n + 2 N) \tr[B]^2)$
		
				\item $\mathbb{E}[\tr[BQQBQQ]] = N((20 + n(11 + n) + 21 N + n(7 + n)N + 2(3 + n)N^2 + N^3)\tr[B^2] + (10 + 5n + n^2 + 5(2 + n)N + 4 N^2) \tr[B]^2)$
		
	\item $\mathbb{E}[\tr[BQ]\tr[BQ]] = (N \tr[B]^2+2 \tr[B^2])N$
	
		\item $\mathbb{E}[\tr[BQ]\tr[BQQ]] = N(2(2 + n + 2N)\tr[B^2] + (2 + N(1 + n + N))\tr[B]^2)$
		
		\item $\mathbb{E}[\tr[BQQ]\tr[BQQ]] =N(2(10 + 5n + n^2 + 5(2 + n)N + 4N^2)\tr[B^2] + (n^2 N + 2n(3 + N + N^2) + (1 + N)(10 + N + N^2))\tr[B]^2)$
		
		The following assume that $q_i $ and $q_j$ are n-dimensional orthonormal vectors uniformly distributed on the surface of a sphere. 
		
	\item $\mathbb{E}[(b^T q_i)^2] = \norm{b}^2 \frac{1}{n}$

  \item $\mathbb{E}[(b^T q_i)^4] = \norm{b}^4 \frac{3}{n(n+2)}$
  
    \item $\mathbb{E}[(b^T q_i)^2 (b^T q_j)^2] = \norm{b}^4 \frac{1}{n(n+2)}$
\end{enumerate}

\subsection{Some Proofs}

\subsubsection{$\mathbb{E}[\tr[BQBQQ]] $}

\begin{align*}
\mathbb{E}[BQBQQ] &= \mathbb{E}\left[B \left(\sum_{i=1}^N Q_i \right)B \left(\sum_{j=1}^N Q_j \right) \left(\sum_{k=1}^N Q_k \right)\right]\\
&= N \mathbb{E}\left[B Q_1 B \left(\sum_{j=1}^N Q_j \right) \left(\sum_{k=1}^N Q_k \right)\right] \\
&= N \left( \mathbb{E}\left[B Q_1 B Q_1 Q_1  + (N-1) B Q_1 B Q_2 Q_2 + 2 (N-1)  B Q_1 B Q_1 Q_2 + (N-1)(N-2) B Q_1 B Q_2 Q_3 \right] \right) \\
\tr[\mathbb{E}[BQBQQ] ] &= N \left((n+4) \tr[B]^2+ (8+2n) \tr[B^2]  + (N-1) \tr[B^2] (n+2) + 2 (N-1)  (2 \tr[B^2] + \tr[B]^2) + (N-1)(N-2) \tr[B^2] \right) \\
&=N((4 + n + (3 + n)N + N^2)\tr[B^2] + (2 + n + 2 N) \tr[B]^2)
\end{align*}

\subsubsection{$\mathbb{E}[\tr[BQQBQQ]] $}

\begin{align*}
\mathbb{E}[BQQBQQ]] &= \mathbb{E}\left[B \left(\sum_{i=1}^N Q_i \right) \left(\sum_{j=1}^N Q_j \right)B \left(\sum_{k=1}^N Q_k \right) \left(\sum_{l=1}^N Q_l \right)\right]\\
&= N \mathbb{E}\left[B \left(\sum_{i=1}^N Q_i \right) \left(\sum_{j=1}^N Q_j \right)B \left(\sum_{k=1}^N Q_k \right)  Q_1 \right]
\end{align*}

This requires quartic products of outer products of Gaussians. To track the permutations, we take N=4 and represent k slices of the i,j interactions. The rows of the slices represent $Q_i$ and columns represent $Q_j$. The color coding distinguishes the eight unique interactions.

\begin{tiny}
\begin{tabular}{|c|c|c|c|}
	\hline
	\cellcolor{red!25} I & III & III & III \\
    \hline
	III & \cellcolor{blue!25}II & \cellcolor{green!25}IV& \cellcolor{green!25}IV \\
	\hline
	III & \cellcolor{green!25}IV & \cellcolor{blue!25} II & \cellcolor{green!25}IV \\
	\hline
	III & \cellcolor{green!25}IV & \cellcolor{green!25}IV & \cellcolor{blue!25}II \\
	\hline
\end{tabular}
\quad
\begin{tabular}{|c|c|c|c|}
	\hline
 III &	\cellcolor{black!25}  V & 	\cellcolor{red!55} VI & 	\cellcolor{red!55} VI \\
	\hline
	\cellcolor{orange!50} 	VII & III &	\cellcolor{red!55} VI & 	\cellcolor{red!55} VI \\
	\hline
		\cellcolor{red!55}	VI &	\cellcolor{red!55} VI &\cellcolor{green!25} IV & \cellcolor{yellow!50}VIII \\
	\hline
	\cellcolor{red!55} 	VI & 	\cellcolor{red!55} VI & \cellcolor{yellow!50}VIII & \cellcolor{green!25}IV \\
	\hline
\end{tabular}
\quad
\begin{tabular}{|c|c|c|c|}
	\hline
	III &	\cellcolor{red!55}  VI & 	\cellcolor{black!25} V & 	\cellcolor{red!55} VI \\
	\hline
	\cellcolor{red!55} 	VI & \cellcolor{green!25} IV &	\cellcolor{red!55} VI & 	\cellcolor{yellow!50} VIII \\
	\hline
	\cellcolor{orange!50}	VII &	\cellcolor{red!55} VI & III & \cellcolor{red!55}VI \\
	\hline
	\cellcolor{red!55} 	VI & 	\cellcolor{yellow!50} VIII & \cellcolor{red!55}VI & \cellcolor{green!25}IV \\
	\hline
\end{tabular}
\quad
\begin{tabular}{|c|c|c|c|}
	\hline
	III &  	\cellcolor{red!55} VI & 	\cellcolor{red!55} VI & 		\cellcolor{black!25} V \\
	\hline
	\cellcolor{red!55} 	VI & \cellcolor{green!25} IV &		\cellcolor{yellow!50} VIII & 	\cellcolor{red!55} VI \\
	\hline
	\cellcolor{red!55} 	VI&	\cellcolor{yellow!50} VIII & \cellcolor{green!25} IV  & \cellcolor{red!55}VI \\
	\hline
	\cellcolor{orange!50}	VII &	\cellcolor{red!55} 	VI  & \cellcolor{red!55}VI &  III \\
	\hline
\end{tabular}
\end{tiny}
\begin{align*}
\tr[\mathbb{E}[BQQBQQ]] &= N \tr\left[\mathbb{E}\left[B \left(\sum_{i=1}^N Q_i \right) \left(\sum_{j=1}^N Q_j \right)B \left(\sum_{k=1}^N Q_k \right) Q_1 \right]\right]\\
&= N \tr[ \mathbb{E}\left[B Q_1 Q_1 B Q_1  Q_1 \right]] \\
&+N (N-1)  \tr[\mathbb{E}\left[B Q_2 Q_2 B Q_1  Q_1 \right]] \\
&+4 N (N-1)  \tr[\mathbb{E}\left[B Q_1 Q_2 B Q_1  Q_1 \right]] \\
&+2 N (N-1) (N-2)  \tr[\mathbb{E}\left[B Q_2 Q_3 B Q_1  Q_1 \right]] \\
&+N (N-1)  \tr[\mathbb{E}\left[B Q_1 Q_2 B Q_2  Q_1 \right]] \\
&+4 N (N-1)(N-2) \tr[\mathbb{E}\left[B Q_1 Q_3 B Q_2  Q_1 \right]] \\
&+N (N-1) \tr[\mathbb{E}\left[B Q_2 Q_1 B Q_2  Q_1 \right]] \\
&+N (N-1)(N-2)(N-3) \tr[\mathbb{E}\left[B Q_4 Q_3 B Q_2  Q_1 \right]] \\
&= N ((n^2+24+10 n) \tr[B]^2+ (2n^2+20n+48) \tr[B^2])\\
&+N (N-1)  \tr[B^2](2+n)^2 \\
&+4 N (N-1)  ((n+4) \tr[B]^2+ (8+2n) \tr[B^2])\\
&+2 N (N-1) (N-2) \tr[B^2](2+n) \\
&+N (N-1)  (4 \tr[B^2]+\tr[B]^2(4+n)) \\
&+4 N (N-1)(N-2)(2 \tr[B^2] + \tr[B]^2) \\
&+N (N-1) (\tr[B^2] (6 + n) + 2 \tr[B]^2 )\\
&+N (N-1)(N-2)(N-3) \tr[B^2].
\end{align*}

\subsubsection{$\mathbb{E}[\tr[BQQ]\tr[BQ]]$}

\begin{align*}
&=N \mathbb{E}[(x_1^T BQ x_1+x_2^T B Q x_2+...+x_N^T BQ x_N)x_1^T B x_1] \\
&= N \mathbb{E}[(x_1^T B (x_1 x_1^T+x_2 x_2^T+...+x_N x_N^T) x_1 x_1^T B x_1]\\
&+ N \mathbb{E}[(x_2^T B (x_1 x_1^T+x_2 x_2^T+...+x_N x_N^T) x_2 x_1^T B x_1]\\
&+ .... \\
&+ N \mathbb{E}[(x_N^T B (x_1 x_1^T+x_2 x_2^T+...+x_N x_N^T) x_N x_1^T B x_1]\\
&= N \mathbb{E}[(x_1^T B (x_1 x_1^T+x_2 x_2^T+...+x_N x_N^T) x_1 x_1^T B x_1]\\
&+ N (N-1) \mathbb{E}[(x_2^T B (x_1 x_1^T+x_2 x_2^T+...+x_N x_N^T) x_2 x_1^T B x_1]\\
&= N \mathbb{E}[(x_1^T B (x_1 x_1^T+x_2 x_2^T+...+x_N x_N^T) x_1 x_1^T B x_1]\\
&+ N (N-1) \mathbb{E}[x_2^T B x_1 x_1^T x_2 x_1^T B x_1] \\
&+ N (N-1) \mathbb{E}[x_2^T B x_2 x_2^T x_2 x_1^T B x_1] \\
&+ N (N-1) (N-2) \mathbb{E}[x_2^T B x_3 x_3^T x_2 x_1^T B x_1] \\
&=N ((n+4) \tr[B]^2+ (8+2n) \tr[B^2] + (N-1)(2\tr[B^2]+\tr[B]^2)) \\
&+ N (N-1) (2 \tr[B^2]+\tr[B]^2)\\
&+ N (N-1) (2+n)\tr[B]^2\\
&+ N (N-1)(N-2) \tr[B]^2
\end{align*}

\subsubsection{$\mathbb{E}[\tr[BQQ]\tr[BQQ]]$}
This involves terms of the form
$
\mathbb{E}[x_j^T B x_i x_i^T x_j x_1^T B x_k x_k^T x_1]
$

Term I
\begin{align*}
\mathbb{E}[x_1^T B x_1 x_1^T x_1 x_1^T B x_1 x_1^T x_1] = (n^2+24+10 n) \tr[B]^2+ (2n^2+20n+48) \tr[B^2]
\end{align*}

Term II
\begin{align*}
\mathbb{E}[x_2^T B x_2 x_2^T x_2 x_1^T B x_1 x_1^T x_1] = \tr[B]^2(n+2)^2
\end{align*}

Term III
\begin{align*}
\mathbb{E}[x_2^T B x_1 x_1^T x_2 x_1^T B x_1 x_1^T x_1] = \mathbb{E}[  x_1^T B x_1 x_1^T B x_1 x_1^T x_1]  = (n+4) \tr[B]^2+ (8+2n) \tr[B^2]
\end{align*}

Term IV
\begin{align*}
\mathbb{E}[x_3^T B x_2 x_2^T x_3 x_1^T B x_1 x_1^T x_1] = \mathbb{E}[x_2^T B x_2  x_1^T B x_1 x_1^T x_1] = \tr[B]^2(2+n)
\end{align*}

Term V
\begin{align*}
\mathbb{E}[x_2^T B x_1 x_1^T x_2 x_1^T B x_2 x_2^T x_1]& =\mathbb{E}[x_2^T C x_2 a^T x_2 x_2^T b] = \mathbb{E}[a^T (C+C^T)b + a^T \tr[C] b] \\
&= \mathbb{E}[x_1^T B (B x_1 x_1^T + x_1 x_1^T B)x_1 + x_1^T B \tr[B x_1 x_1^T]x_1]\\
&= \mathbb{E}[x_1^T B^2  x_1 x_1^T x_1 + 2 x_1^T B x_1  x_1^T B x_1]\\
&= \tr[B^2](2+n) + 4 \tr[B^2] + 2 \tr[B]^2
\end{align*}

Term VI
\begin{align*}
\mathbb{E}[x_3^T B x_1 x_1^T x_3 x_1^T B x_2 x_2^T x_1] = \mathbb{E}[x_3^T B x_1 x_1^T x_3 x_1^T B  x_1]  = \mathbb{E}[x_1^T B x_1 x_1^T B  x_1] = 2 \tr[B^2] + \tr[B]^2  
\end{align*}

Term VII
\begin{align*}
\mathbb{E}[x_1^T B x_2 x_2^T x_1 x_1^T B x_2 x_2^T x_1]& = \mathbb{E}[\tr[B x_2 x_2^T (B x_2 x_2^T+ x_2 x_2^T B)] + \tr[B x_2 x_2^T]^2 ] \\
& = \mathbb{E}[\tr[B x_2 x_2^T (B x_2 x_2^T+ x_2 x_2^T B)] + \tr[B x_2 x_2^T]^2 ] \\
& = 2 \mathbb{E}[x_2^T B x_2 x_2^T B x_2 ] + \mathbb{E}[x_2^T x_2 x_2^T B^2 x_2] \\
& = 4 \tr[B^2] + 2 \tr[B]^2 + (2+n)\tr[B^2] \\
\end{align*}

Term VIII
\begin{align*}
\mathbb{E}[x_4^T B x_3 x_3^T x_4 x_1^T B x_2 x_2^T x_1] &= \mathbb{E}[x_4^T B  x_4 x_1^T B x_1]  = \tr[B]^2 
\end{align*}

\subsubsection{$\mathbb{E}[\tr[Q]b^TQQb]] $}

\begin{align*}
 &= N \mathbb{E}\left[ x_1^T x_1 b^T \left(\sum_{j=1}^N Q_j \right) \left(\sum_{k=1}^N Q_k \right)b\right] \\
&= N \left( \mathbb{E}\left[x_1^T x_1 b^T x_1 x_1^T x_1 x_1^T b  + (N-1) x_1^T x_1 b^T x_2 x_2^T x_2 x_2^T b + 2 (N-1)  x_1^T x_1 b^T x_1 x_1^T x_2 x_2^T b + (N-1)(N-2) x_1^T x_1 b^T x_2 x_2^T x_3 x_3^T b \right] \right) \\
&= N b^T b \left( (n^2+6n+8)+ (N-1)(n(2+n)) + 2 (N-1)(2+n)   + (N-1)(N-2)n  \right) \\
&= b^T b N (1 + n + N) (4 + n N)
\end{align*}

\bibliographystyle{siamplain}
\bibliography{refs}

\end{document}